\documentclass{amsart}
\usepackage{amssymb, amsthm, amsmath}

\newtheorem{theorem}{Theorem}[section]

\newtheorem{prop}[theorem]{Proposition}

\theoremstyle{definition}
\newtheorem{remark}[theorem]{Remark}
\newtheorem{question}[theorem]{Question}
\newtheorem{conj}[theorem]{Conjecture}

\def\Hom{\mathop{\mathrm{Hom}}}
\def\Ext{\mathop{\mathrm{Ext}}}
\def\ann{\mathop{\mathrm{ann}}}
\def\Ass{\mathop{\mathrm{Ass}}}
\def\Spec{\mathop{\mathrm{Spec}}}
\def\Supp{\mathop{\mathrm{Supp}}}
\def\Ker{\mathop{\mathrm{Ker}}}
\def\Im{\mathop{\mathrm{Im}}}
\def\m{\mathfrak{m}}
\def\p{\mathfrak{p}}
\def\P{\mathfrak{P}}
\def\a{\mathfrak{a}}
\def\N{\mathbb{N}}
\def\C{\mathbb{C}}
\def\Z{\mathbb{Z}}
\def\Pdot{P_{\bullet}}

\begin{document}

\title{Associated primes of local cohomology modules}

\author{Anurag K. Singh}

\address{School of Mathematics \\
686 Cherry Street \\
Georgia Institute of Technology \\
Atlanta, GA 30332-0160, USA.
E-mail: {\tt singh@math.gatech.edu}}

\thanks {The author was supported in part by grants from the National
Science Foundation.}

\date{\today}

\maketitle

\section{introduction}

Throughout, $R$ will denote a commutative Noetherian ring with a unit element.
Let $\a$ be an ideal of $R$, and $i$ a non-negative integer. The {\it local
cohomology module}\/ $H^i_\a(R)$ is defined as
$$
H^i_\a(R) = \varinjlim_{k \in \N} {\Ext}^i_R \left( R/\a^k,R \right),
$$
where the maps in the direct limit system are those induced by the natural
surjections $R/\a^{k+1} \longrightarrow R/\a^k$. If $\a$ is generated by
elements $x_1,\dots,x_n,$ then $H^i_\a(R)$ is isomorphic to the $i$\/th
cohomology module of the extended \v Cech complex
$$
0 \longrightarrow R \longrightarrow \bigoplus_{i=1}^n R_{x_i} \longrightarrow 
\bigoplus_{i<j} R_{x_i x_j} \longrightarrow \cdots \longrightarrow
R_{x_1 \cdots x_n} \longrightarrow 0.
$$
For an element $f \in R$ and a positive integer $m$, we use 
$[f+(x_1^m, \dots, x_n^m)]$ to denote the cohomology class
$$
\left[\frac{f}{x_1^m \cdots x_n^m} \right] \ \in \ 
\frac{R_{x_1 \cdots x_n}}{\sum_i R_{x_1 \cdots \hat{x}_i \cdots x_n}}
\ \cong \ H^n_\a(R).
$$
It is easily seen that $[f+(x_1^m, \dots, x_n^m)]=0$ in $H^n_\a(R)$ if and only
if there exists an integer $k \ge 0$, such that
$$
fx_1^k\cdots x_n^k \in \big(x_1^{m+k},\dots,x_n^{m+k}\big)R.
$$
Consequently $H^n_\a(R)$ may be also identified with the direct limit
$$
\varinjlim_{m \in \N} R/(x_1^m,\dots,x_n^m)R,
$$
where the map $R/(x_1^m,\dots,x_n^m) \longrightarrow
R/(x_1^{m+1},\dots,x_n^{m+1})$ is multiplication by the image of the element
$x_1 \cdots x_n$.

As these descriptions suggest, $H^i_\a(R)$ is usually not finitely generated as
an $R$-module. However local cohomology modules have useful finiteness
properties in certain cases, e.g., for a local ring $(R,\m)$, the modules
$H^i_\m(R)$ satisfy the descending chain condition. This implies, in particular,
that for all $i \ge 0$,
$$
{\Hom}_R \left( R/\m, H^i_\m(R) \right) \cong 0 :_{H^i_\m(R)} \m
$$
is a finitely generated $R$-module. Grothendieck conjectured that for all
ideals $\a \subset R$, the modules
$$
{\Hom}_R \left( R/\a, H^i_\a(R) \right) \cong 0 :_{H^i_\a(R)} \a 
$$
are finitely generated, \cite[Expos\'e XIII, page 173]{SGA2}. In \cite[\S 3]{Ha}
Hartshorne gave a counterexample to this conjecture: Let $K$ be a field and $R$
be the hypersurface 
$$
K[w,x,y,z]/(wx-yz).
$$
Set $\a=(x,y)$ and consider the local cohomology module $H^2_\a(R)$. It is
easily seen that the elements
$$
[y^nz^n+(x^{n+1},y^{n+1})R] \in H^2_\a(R) \quad \text{for} \quad n \ge 0 
$$
are nonzero, and are killed by the maximal ideal $\m=(w,x,y,z)$. In fact, they
span the module $0 :_{H^2_\a(R)} \m$ which is a vector space of countably
infinite dimension, and so it cannot be finitely generated as an $R$-module. It
follows that $0 :_{H^2_\a(R)} \a$ is not finitely generated as well.

In \cite{Ha} Hartshorne also began the study of the cofiniteness of local
cohomology modules: An $R$-module $M$ is $\a$-{\it cofinite} if
$\Supp(M) \subseteq V(\a)$ and ${\Ext}^i_R(R/\a,M)$ is finitely generated for
all $i \ge 0$. Some of the work on cofiniteness may be found in the papers
\cite{Ch, DM, HK, HM, Kw, Me, Ya}, and \cite{Yo}. A related question on the
torsion in local cohomology modules was raised by Huneke at the Sundance
Conference in 1990, and will be our main focus here.

\begin{question}\cite{Hu1}\label{Huneke}
Is the number of associated prime ideals of a local cohomology module
$H_\a^i(R)$ always finite?
\end{question}

The first results were obtained by Huneke and Sharp.

\begin{theorem}\cite[Corollary 2.3]{HS}
Let $R$ be a regular ring containing a field of positive characteristic, and
$\a \subset R$ an ideal. Then for all $i \ge 0$,
$$
\Ass H^i_\a(R) \subseteq \Ass {\Ext}^i_R(R/\a,R) \qquad \qquad (*)
$$
In particular, $\Ass H^i_\a(R)$ is a finite set.
\end{theorem}

\begin{remark}
The proof of the above theorem relies heavily on the flatness of the Frobenius
endomorphism which, by \cite[Theorem 2.1]{Ku}, characterizes regular rings of
positive characteristic. The containment $(*)$ may fail for regular rings of
characteristic zero: Let $R=\C[u,v,w,x,y,z]$, and $\a$ be the ideal generated by
the $2 \times 2$ minors $\Delta_i$ of the matrix
$$
M = \begin{pmatrix}
u & v & w \\
x & y & z
\end{pmatrix}.
$$
Then ${\Ext}^3_R(R/\a,R) = 0$ since $R/\a$ has projective dimension two as an
$R$-module. However, as observed by Hochster, the module $H^3_\a(R)$ is nonzero:
To see this, consider the linear action of $G=SL_2(\C)$ on $R$, where an element
$g \in G$ maps the entries of the matrix $M$ to those of the matrix
$g \times M$. The ring of invariants for this action is the polynomial ring
$R^G =\C[\Delta_1, \Delta_2, \Delta_3]$. Since $SL_2(\C)$ is linearly reductive,
the inclusion $R^G \hookrightarrow R$ splits via an $R^G$-linear retraction,
and so
$$
H^3_{(\Delta_1, \Delta_2, \Delta_3)} (R^G) \longrightarrow H^3_\a (R)
$$
is a split inclusion. Since the module
$H^3_{(\Delta_1, \Delta_2, \Delta_3)}(R^G)$ is nonzero, it follows that
$H^3_\a (R)$ must be nonzero as well.
\end{remark}

While $\Ass H^i_\a(R)$ may not be a subset of $\Ass {\Ext}^i_R(R/\a,R)$,
Question \ref{Huneke} does have an affirmative answer for all unramified
regular local rings by combining the result of Huneke-Sharp with the following
two theorems of Lyubeznik.

\begin{theorem}\cite[Corollary 3.6 (c)]{Ly1}
Let $R$ be a regular ring containing a field of characteristic zero and $\a$ an
ideal of $R$. Then for every maximal ideal $\m$ of $R$, the set of associated
primes of a local cohomology module $H^i_\a(R)$, which are contained in the
ideal $\m$, is finite.

If the regular ring $R$ is finitely generated over a field of characteristic
zero, then $\Ass H^i_\a(R)$ is a finite set.
\end{theorem}

To illustrate the key point here, consider the case where $R=\C[x_1,\dots,x_n]$,
and let $D$ be the ring of $\C$-linear differential operators on $R$. It turns
out that $D$ is left and right Noetherian, that $H^i_\a(R)$ is a finitely
generated $D$-module, and consequently that $\Ass H^i_\a(R)$ is finite.
Lyubeznik's result below also uses $D$-modules, though the situation in mixed
characteristic is more subtle.

\begin{theorem}\cite[Theorem 1]{Ly2}
If $R$ is an unramified regular local ring of mixed characteristic, and $\a$ is
an ideal of $R$, then $\Ass H^i_\a(R)$ is a finite set.
\end{theorem}

So far we have restricted the discussion to local cohomology modules of the form
$H^i_\a(R)$. For an $R$-module $M$, the local cohomology modules $H^i_\a(M)$ are
defined similarly as
$$
H^i_\a(M) = \varinjlim_{k \in \N} {\Ext}^i_R \left( R/\a^k,M \right), \qquad
\text{where} \qquad i \ge 0.
$$
If $M$ is a finitely generated $R$-module, then $H^0_\a(M)$ may be identified
with the submodule of $M$ consisting of elements which are killed by a power of
the ideal $\a$, and consequently $H^0_\a(M)$ is a finitely generated $R$-module.
If $i$ is the smallest integer for which $H_\a^i(M)$ is not finitely generated,
then the set $\Ass H_\a^i(M)$ is also finite, as proved in \cite{BF} and
\cite{KS}. Other positive answers to Question \ref{Huneke} include results in
small dimensions such as the following theorem due to Marley:

\begin{theorem}\cite[Corollary 2.7]{Ma}
Let $R$ be a local ring and $M$ a finitely generated $R$-module of dimension at
most three. Then $\Ass H_\a^i(M)$ is finite for all ideals $\a \subset R$.
\end{theorem}

For some of the other work on this question, we refer the reader to the papers
\cite{BKS, BRS, He, Ly3} and \cite{MV}.

\section{$p$-torsion}

In \cite{ptor1} the author constructed a hypersurface for which a local
cohomology module has infinitely many associated prime ideals, thereby
demonstrating that Question \ref{Huneke}, in general, has a negative answer.
Since the argument is quite elementary, we include it here.

\begin{theorem}\cite[\S 4]{ptor1}
Consider the hypersurface
$$
R=\Z[u,v,w,x,y,z]/(ux+vy+wz)
$$
and the ideal $\a=(x,y,z)R$. Then for every prime integer $p$, the local
cohomology module $H_\a^3(R)$ has a $p$-torsion element. Consequently
$H_\a^3(R)$ has infinitely many associated prime ideals.
\end{theorem}

\begin{proof}
We identify $H_\a^3(R)$ with the direct limit
$$
\varinjlim_{k \in \N} R/(x^k, y^k, z^k)R,
$$
where the maps are induced by multiplication by the element $xyz$. For a prime
integer $p$, the fraction
$$
\lambda_p=\frac{(ux)^p + (vy)^p + (wz)^p}{p}
$$
has integer coefficients, and is therefore an element of $R$. We claim that the
element
$$
\eta_p = [\lambda_p + (x^p, y^p, z^p)R] \in H_\a^3(R)
$$
is nonzero and $p$-torsion. Note that $p \cdot \eta_p =
[p\lambda_p + (x^p, y^p, z^p)R] = 0$, and what remains to be checked is that
$\eta_p$ is nonzero, i.e., that
$$
\lambda_p (xyz)^k \notin \left( x^{p+k}, y^{p+k}, z^{p+k} \right)R \quad \text
{for all} \quad k \in \N. 
$$
We assign weights to the $\Z$-algebra generators of the ring $R$ as follows:
\begin{align*}
x: (1,0,0,0), \qquad\qquad\qquad & u: (-1,0,0,1), \\
y: (0,1,0,0), \qquad\qquad\qquad & v: (0,-1,0,1), \\
z: (0,0,1,0), \qquad\qquad\qquad & w: (0,0,-1,1). 
\end{align*}
With this grading, $\lambda_p$ is a homogeneous element of degree $(0,0,0,p)$. 
Now suppose we have a homogeneous equation of the form
$$
\lambda (xyz)^{k}=c_1x^{p+k}+c_2y^{p+k}+c_3z^{p+k},
$$
then we must have $\deg(c_1) = (-p,k,k,p)$, i.e., $c_1$ must be an integer
multiple of the monomial $u^p y^k z^k$. Similarly $c_2$ is an integer multiple
of $v^p z^k x^k$ and $c_3$ of $w^p x^k y^k$. Consequently
\begin{align*}
\lambda (xyz)^k 
& \in \left( u^p y^k z^k x^{p+k}, \ v^p z^k x^k y^{p+k},\ w^p x^k y^k z^{p+k}
\right) R \\
& = (xyz)^k \left( u^px^p, \ v^py^p, \ w^pz^p \right) R,
\end{align*}
and so $\lambda \in (u^px^p, v^py^p, w^pz^p)R$. After specializing
$u \mapsto 1, v \mapsto 1, w \mapsto 1$, this implies that
$$
\frac{x^p+y^p+(-1)^p(x+y)^p}{p} \in \left( p, x^p, y^p \right) \Z[x,y],
$$
which is easily seen to be false.
\end{proof}

This example, however, does not shed light on Question \ref{Huneke} in the case
of local rings or rings containing a field. Katzman constructed the first
examples to demonstrate that Huneke's question has a negative answer in these
cases as well, \cite{Ka2}. The equicharacteristic case is discussed here in 
\S \ref{equicharacteristic}. We next recall a conjecture of Lyubeznik.

\begin{conj}\cite[Remark 3.7 (iii)]{Ly1}
If $R$ is a regular ring and $\a$ an ideal, then the local cohomology modules
$H_\a^i(R)$ have finitely many associated prime ideals.
\end{conj}

This conjecture has been settled for unramified regular local rings by the
results of Huneke-Sharp and Lyubeznik mentioned earlier. However it remains open
for polynomial rings over the integers, and we discuss some of its implications
in this case.

\begin{remark}
Let $R$ be a polynomial ring in finitely many variables over the integers, and
let $\a$ be an ideal of $R$. Then for every prime integer $p$, we have a short
exact sequence
$$
0 \longrightarrow R \overset{p}\longrightarrow R \longrightarrow R/pR
\longrightarrow 0,
$$
which induces a long exact sequence of local cohomology modules,
$$
\cdots \longrightarrow H_\a^{i-1} (R/pR) \overset{\delta_p^{i-1}}
\longrightarrow H_\a^i (R) \overset{p} \longrightarrow H_\a^i (R)
\longrightarrow H_\a^i (R/pR) \overset{\delta_p^i} \longrightarrow
H_\a^{i+1} (R) \overset{p} \longrightarrow \cdots .
$$
The image of each connecting homomorphism $\delta_p^i$ is annihilated by $p$,
and hence every nonzero element of $\delta_p^i(H_\a^i(R/pR))$ is a $p$-torsion
element. Consequently Lyubeznik's conjectures implies that for all but finitely
many prime integers $p$, we must have $\delta_p^{i} = 0$ for all $i \ge 0$.
\end{remark}

\begin{remark}
Again, let $R$ be a polynomial ring over the integers. Let $f_i,g_i$ be elements
of $R$ such that
$$
f_1g_1 + f_2g_2 + \dots + f_ng_n = 0.
$$
Consider the ideal $\a = (g_1, \dots, g_n)R$ and the local cohomology module
$$
H_\a^n(R) = \varinjlim_{k \in \N} R/(g_1^k, \dots, g_n^k)R,
$$
where the maps in the direct system are induced by multiplication by the element
$g_1 \cdots g_n$. For a prime integer $p$ and prime power $q=p^e$, let
$$
\lambda_q = \frac{(f_1g_1)^q + \dots + (f_ng_n)^q}{p}.
$$
Then $\lambda_q \in R$, and we set
$$
\eta_q = [\lambda_q + (g_1^q, \dots, g_n^q)R] \in H_\a^n(R).
$$
It is immediately seen that $p \cdot \eta_q=0$ and so if $\eta_ q$ is a nonzero 
element of $H_\a^n(R)$, then it must be a $p$-torsion element. Hence
Lyubeznik's conjecture implies that for all but finitely many prime integers
$p$, the elements $\eta_q$ must be zero, i.e., for some $k \in \N$, which may
depend on $q=p^e$, we have
$$
\lambda_q (g_1 \cdots g_n)^k \in (g_1^{q+k}, \dots, g_n^{q+k})R.
$$
\end{remark}

This motivates the following conjecture:
\begin{conj}
Let $R$ be a polynomial ring over the integers, and let $f_i,g_i$ be elements of
$R$ such that
$$
f_1g_1 + \dots + f_ng_n = 0.
$$
Then for every prime power $q=p^e$, there exists $k \in \N$ such that 
$$
\frac{(f_1g_1)^q + \dots + (f_ng_n)^q}{p} (g_1 \cdots g_n)^{k} \in
(g_1^{q+k}, \dots, g_n^{q+k})R.
$$
\label{conj}
\end{conj}

The above conjecture is easily established if $n=2$, or if the elements
$g_1, \dots, g_n$ form a regular sequence. The conjecture is also true if
$n=3$, provided the elements $f_1, f_2, f_3$ form a regular sequence:

\begin{theorem}\cite[Theorem 2.1]{ptor2}
Let $R$ be a polynomial ring over the integers and $f_i, g_i$ be elements of
$R$ such that $f_1, f_2, f_3$ form a regular sequence in $R$ and
$$
f_1g_1 + f_2g_2 + f_3g_3 = 0.
$$
Let $q=p^e$ be a prime power. Then for $k=q-1$, we have
$$
\frac{(f_1g_1)^q + (f_2g_2)^q + (f_3g_3)^q}{p} (g_1 g_2 g_3)^k \in
(g_1^{q+k}, g_2^{q+k}, g_3^{q+k})R.
$$
\end{theorem}

\section{The equicharacteristic case}\label{equicharacteristic}

Recently Katzman constructed the following example in \cite{Ka2}: Let $K$ be an
arbitrary field, and consider the hypersurface
$$
R = K[s,t,u,v,x,y]/\big(su^2x^2-(s+t)uxvy+tv^2y^2\big).
$$
Katzman showed that the local cohomology module $H^2_{(x,y)}(R)$ has infinitely
many associated prime ideals. Since the defining equation of this hypersurface
factors as
$$
su^2x^2-(s+t)uxvy+tv^2y^2 = (sux-tvy)(ux-vy),
$$
the ring in Katzman's example is not an integral domain. In \cite{SS} Swanson
and the author generalize Katzman's construction and obtain families of
examples which include examples over normal domains and, in fact, over
hypersurfaces with rational singularities:

\begin{theorem}\cite[Theorem 1.1]{SS}\label{main}
Let $K$ be an arbitrary field, and consider the hypersurface
$$
S=\frac{K[s,t,u,v,w,x,y,z]}{\big(su^2x^2+sv^2y^2+tuxvy+tw^2z^2\big)}
$$
Then $S$ is a standard $\N$-graded normal domain for which the local cohomology
module $H^3_{(x,y,z)}(S)$ has infinitely many associated prime ideals.

If $\m$ denotes the homogeneous maximal ideal $(s,t,u,v,w,x,y,z)$, then 
the local cohomology module $H^3_{(x,y,z)}(S_\m)$ has infinitely many
associated prime ideals as well.

If $K$ has characteristic zero, then $S$ has rational singularities. If $K$ has
positive characteristic, then $S$ is F-regular.
\end{theorem}

While we refer to \cite{SS} for details as well as more general constructions,
we would like to sketch a proof of the above theorem.

\begin{proof}[Sketch of the proof]
Let $K$ be arbitrary field, and consider the hypersurface
$$
A=K[s,t,a,b]/(sa^2+tab+b^2).
$$
We use the $\N$-grading on $A$ where $A_0=K[s,t]$, and $a$ and $b$ have degree
$1$. For an integer $n \ge 2$, the presentation matrix of the $A_0$-module
$[A/(a^n,b^n)]_n$ is the $(n-1) \times (n-1)$ matrix
$$
M_{n-1}=\left[\begin{matrix}
t & s \cr
s & t & s \cr
& \ddots & \ddots & \ddots \cr
& & s & t & s \cr
& & & s & t \cr
\end{matrix}\right],
$$
where the elements $s,t,s$ occur along the three central diagonals, and the
other entries are zero. These are special cases of Toeplitz matrices, and the
determinants of these matrices will be the source of the infinitely many
associated primes.

It is convenient to define $\det M_0=1$, and setting $Q_n(s,t) = \det M_n$, it
is easily seen that we have
$$
Q_0=1, \quad Q_1=t, \quad \text{and} \quad Q_{n+2}=tQ_{n+1}-s^2Q_n \quad
\text{for all} \quad n \ge 0.
$$
This recursion enables us to compute a generating function for the polynomials
$Q_n(s,t)$, and to check that the polynomials $\{Q_n(s,t)\}_{n \in \N}$ have
infinitely many distinct irreducible factors. For example, if $K=\C$ it may be
verified that $Q_n(s,t)$ factors as
$$
Q_n(s,t) = \prod_{r=1}^n \left(t-se^{r\pi i/(n+1)} -se^{-r\pi i/(n+1)}\right)
\quad \text{for all} \quad n \in \N.
$$
This implies that the set
$$
\bigcup_{n \in \N} {\Ass}_{A_0} [A/(a^n,b^n)]_n
$$
is infinite and, in fact, a straightforward computation shows that 
$$
(a^n,b^n):_{A_0} sab^{n-1} = (Q_{n-1})A_0 \quad \text{for all} \quad n \in \N.
$$

We next consider the ring
$$
B=K[s,t,a,b,c]/\big(sa^2+sb^2+tab+tc^2\big)
$$
with the $\N$-grading where $B_0=K[s,t]=A_0$, and $a$, $b$, and $c$ have degree
$1$. Note that $B/cB \cong A$, and so
$$
(a^n,b^n,c):_{B_0} sab^{n-1} = (Q_{n-1})B_0 \quad \text{for all} \quad n \in \N.
$$
We identify $B$ with the subring of $S$ which is generated, as a $K$-algebra,
by the elements $s,t,a=ux,b=vy$, and $c=wz$. For a fixed integer $n \ge 1$, let
$$
\eta_n = \left[s(ux)(vy)^{n-1}z^{n-1}+(x^n,y^n,z^n)\right] \in H^3_{(x,y,z)}(S),
$$
where we are using the identification 
$$
H^3_{(x,y,z)}(S)=\varinjlim_{n \in \N} S/(x^n,y^n,z^n)S.
$$
By a multigrading argument, it may be verified that
$$
{\ann}_{S_0} \eta_n = (a^n,b^n,c)B :_{B_0} sab^{n-1} = (Q_{n-1})B_0
$$
where $S_0 = B_0 = K[s,t]$. Since the polynomials $\{Q_n(s,t)\}_{n \in \N}$
have infinitely many distinct irreducible factors, it follows that the set
$$
{\Ass}_{S_0} H^3_{(x,y,z)}(S)
$$
is infinite. For every prime ideal $\p$ of $S_0$ with
$\p \in {\Ass}_{S_0} H^3_{(x,y,z)}(S)$, there exists a prime ideal
$\P \in \Spec S$ such that $\P \in {\Ass}_S H^3_{(x,y,z)}(S)$ and
$\P \cap S_0 = \p$. Consequently the set $\Ass_S H^3_{(x,y,z)}(S)$ must be
infinite as well.

It remains to verify that the hypersurface $S$ has rational singularities (in
characteristic zero) or is F-regular (in positive characteristic). In \cite{SS}
we show that $S$ is F-regular for an arbitrary field $K$ of positive
characteristic. This implies that for all prime integers $p$, the fiber over
$p\Z$ of the map
$$
\Z \longrightarrow
\frac{\Z[s,t,u,v,w,x,y,z]}{\big(su^2x^2+sv^2y^2+tuxvy+tw^2z^2\big)}
$$
is an F-rational ring. By \cite[Theorem~4.3]{sm-ratsing}, it then follows that
$S$ has rational singularities when $K$ has characteristic zero.

We would like to include here a different proof that $S$ has rational
singularities in characteristic zero based on a result from \cite{SW}. We first
note that
$$
S \cong B[u,v,w,x,y,z]/(ux-a,vy-b,wz-c),
$$
and that $B$ is a normal domain. By a result of \cite{BS}, if a local
(or graded) domain $R$ has rational singularities, then so does $R[u,x]/(ux-a)$,
where $a \neq 0$ is a (homogeneous) element of $R$, see also
\cite[Lemma 3.3]{HWY}. By repeated use of this, to show that $S$ has rational
singularities, it suffices to show that the subring $B$ has rational
singularities. In \cite{SW} we obtain a criterion for multigraded rings to have
rational singularities. The bigraded case of this criterion is:

\begin{theorem}
Let $R$ be a normal $\N^2$-graded ring where $R_{\bf 0}$ is a field of
characteristic zero, and $R$ is generated over $R_{\bf 0}$ by elements of
degrees $(1,0)$ and $(0,1)$. Then $R$ has rational singularities if and only if
\begin{itemize}
\item[(i)] $R$ is a Cohen-Macaulay ring for which the multigraded
${\bf a}$-invariant satisfies ${\bf a}(R) < {\bf 0}$, and
\item[(ii)] the localizations $R_\p$ have rational singularities for all
primes $\p$ in the set
$$
\Spec R \setminus V(R_{++}), \qquad \text{where} \qquad 
R_{++} = \bigoplus_{i>0,j>0}R_{i,j}. 
$$
\end{itemize}
\end{theorem}

To apply the theorem, we consider the $\N^2$-grading on $B$ where $s$ and $t$
have degree $(1,0)$ and $a$, $b$, and $c$ have degree $(0,1)$. Then 
${\bf a}(B)=(-1,-1)$, and a straightforward computation using the Jacobian
criterion shows that $B_\p$ is regular for all primes
$\p \in \Spec B \setminus V(B_{++})$.
\end{proof}

\section{An application} \label{tc}

Let $R$ be a ring of characteristic $p>0$, and $R^\circ$ denote the complement
of the minimal primes of $R$. For an ideal $\a=(x_1,\dots,x_n)$
of $R$ and a prime power $q=p^e$, we use the notation
$\a^{[q]}=(x_1^q,\dots,x_n^q)$. The {\it tight closure}\/ of $\a$ is the ideal
$$
\a^* = \{z \in R : \text{ there exists } c\in R^\circ \text{ for which }
cz^q \in \a^{[q]} \text{ for all } q\gg 0 \},
$$
see \cite{HHjams}. A ring $R$ is {\it F-regular}\/ if $\a^* = \a$ for all
ideals $\a$ of $R$ and its localizations. 

More generally, let $F$ denote the Frobenius functor, and $F^e$ its $e$\/th
iteration. If an $R$-module $M$ has presentation matrix $(a_{ij})$, then
$F^e(M)$ has presentation matrix $(a_{ij}^q)$, where $q=p^e$. For modules
$N \subseteq M$, we use $N^{[q]}_M$ to denote the image of $F^e(N) \to F^e(M)$. 
We say that an element $m \in M$ is in the {\it tight closure of $N$ in $M$},
denoted $N^*_M$, if there exists an element $c \in R^\circ$ such that
$cF^e(m) \in N^{[q]}_M$ for all $q \gg 0$. While the theory has found several
applications, the question whether tight closure commutes with localization
remains open even for finitely generated algebras over fields of positive
characteristic.

Let $W$ be a multiplicative system in $R$, and $N \subseteq M$ be finitely
generated $R$-modules. Then 
$$
W^{-1}(N^*_M) \subseteq (W^{-1}N)^*_{W^{-1}M},
$$
where $W^{-1}(N^*_M)$ is identified with its image in $W^{-1}M$. When this
inclusion is an equality, we say that {\it tight closure commutes with
localization at $W$ for the pair $N \subseteq M$}. It may be checked that
this occurs if and only if tight closure commutes with localization at $W$
for the pair $0 \subseteq M/N$. Following \cite{AHH}, we set
$$
G^e(M/N) = F^e(M/N)/(0^*_{F^e(M/N)}).
$$
An element $c \in R^\circ$ is a {\it weak test element}\/ if there exists
$q_0=p^{e_0}$ such that for every pair of finitely generated modules
$N \subseteq M$, an element $m \in M$ is in $N^*_M$ if and only if
$cF^e(m) \in N^{[q]}_M$ for all $q \ge q_0$. By \cite[Theorem 6.1]{HHbasec}, if
$R$ is of finite type over an excellent local ring, then $R$ has a weak test
element. 

\begin{prop}\cite[Lemma 3.5]{AHH} Let $R$ be a ring of characteristic $p>0$ and
$N \subseteq M$ be finitely generated $R$-modules. Then the tight closure of
$N \subseteq M$ commutes with localization at any multiplicative system $W$
which is disjoint from the set $\bigcup_{e\in \N} \Ass F^e(M)/N^{[q]}_M$.

If $R$ has a weak test element, then the tight closure of $N \subseteq M$ also
commutes with localization at multiplicative systems $W$ disjoint from the set
$\bigcup_{e\in \N} \Ass G^e(M/N)$.
\end{prop}

Consider a bounded complex $\Pdot$ of finitely generated projective
$R$-modules,
$$
0 \longrightarrow P_n \overset{d_n}\longrightarrow P_{n-1} \longrightarrow
\cdots \overset{d_1}\longrightarrow P_0 \longrightarrow 0.
$$
The complex $\Pdot$ is said to have {\it phantom homology}\/ at the $i$\/th
spot if
$$
\Ker d_i \subseteq (\Im d_{i+1})^*_{P_i}.
$$
The complex $\Pdot$ is {\it stably phantom acyclic}\/ if $F^e(\Pdot)$ has
phantom homology at the $i$\/th spot for all $i\ge 1$, for all $e \ge 1$. An
$R$-module $M$ has {\it finite phantom projective dimension}\/ if there exits a
bounded stably phantom acyclic complex $\Pdot$ of projective $R$-modules, with
$H_0(\Pdot) \cong M$.

\begin{theorem}\cite[Theorem 8.1]{AHH}
Let $R$ be an equidimensional ring of positive characteristic, which is of
finite type over an excellent local ring. If $N \subseteq M$ are finitely
generated $R$-modules such that $M/N$ has finite phantom projective dimension,
then the tight closure of $N$ in $M$ commutes with localization at $W$ for
every multiplicative system $W$ of $R$.
\end{theorem}

The key points of the proof are that for $M/N$ of finite phantom projective
dimension, the set $\bigcup_e \Ass G^e(M/N)$ has finitely many maximal
elements, and that if $(R,\m)$ is a local ring, then there a positive
integer $B$ such that for all $q=p^e$, the ideal $\m^{Bq}$ kills the local
cohomology module
$$
H^0_\m \left(G^e(M/N) \right).
$$
For more details on this approach to the localization problem, we refer the
reader to the papers \cite{AHH, Ho, Ka1, SN}, and \cite[\S12]{Hu2}.
Specializing to the case where $M=R$ and $N=\a$ is an ideal, we note that
$$
G^e(R/\a) \cong R/ \big(\a^{[q]}\big)^*, \quad \text{where} \quad q=p^e.
$$
Consider the questions:

\begin{question}\label{frob}
Let $R$ be a Noetherian ring of characteristic $p>0$, and $\a$ an ideal of $R$.
\begin{enumerate}
\item Is the set $\bigcup_{q=p^e} \Ass R/ \a^{[q]}$ finite?
\item Is the set $\bigcup_{q=p^e} \Ass R/ \big(\a^{[q]}\big)^*$ finite?
\item For a local domain $(R,\m)$ and an ideal $\a \subset R$, is there a
positive integer $B$ such that
$$
\m^{Bq} H^0_\m \left(R/ (\a^{[q]})^* \right) = 0 \quad \text{for all} \quad
q=p^e \,?
$$
\end{enumerate}
\end{question}

Katzman proved that affirmative answers to Questions \ref{frob}\,$(2)$ and
\ref{frob}\,$(3)$ imply that tight closure commutes with localization:

\begin{theorem}\cite{Ka1}
Assume that for every local ring $(R,\m)$ of characteristic $p>0$ and ideal
$\a \subset R$, the set $\bigcup_q \Ass R/\big(\a^{[q]}\big)^*$ has finitely
many maximal elements. Also, if for every ideal $\a \subset R$, there exists a
positive integer $B$ such that $\m^{Bq}$ kills 
$$
H^0_\m \left(R/ (\a^{[q]})^* \right) \quad \text{for all} \quad q=p^e,
$$
then tight closure commutes with localization for all ideals in Noetherian
rings of characteristic $p>0$.
\end{theorem}

These issues are the source of our interest in associated primes of Frobenius
powers of ideals. It should be mentioned that the situation for {\it ordinary}\/
powers is well-understood: the set $\bigcup_{n \in \N} \Ass R/\a^n$ is
finite for any Noetherian ring $R$, see \cite{Br} or \cite{Ra}. In \cite{Ka1}
Katzman constructed the first example where $\bigcup_{q=p^e} \Ass R/\a^{[q]}$ is
not finite, thereby settling Question \ref{frob}\,$(1)$: For
$$
R = K[t,x,y]/\big(xy(x-y)(x-ty)\big),
$$
he proved that the set $\bigcup_q \Ass R/ (x^q,y^q)$ is infinite. In this
example, and some others, $(x^q,y^q)^*=(x,y)^q$ for all $q=p^e$, and so
$\bigcup_q \Ass R/ (x^q,y^q)^*$ is finite. However we show that Question
\ref{frob}\,$(2)$ also has a negative answer using the local cohomology
examples recorded earlier.

\begin{theorem}[Singh-Swanson]
Let $K$ be a field of characteristic $p>0$, and consider the hypersurface
$$
S=\frac{K[s,t,u,v,w,x,y,z]}{\big(su^2x^2+sv^2y^2+tuxvy+tw^2z^2\big)}.
$$
Then $S$ is F-regular, and the set
$$
\bigcup_{q=p^e} \Ass S/ \big(x^q,y^q,z^q\big) = 
\bigcup_{q=p^e} \Ass S/ \big(x^q,y^q,z^q\big)^*
$$
is infinite.
\end{theorem}

\begin{proof}
The direct system $\{ S/(x^q,y^q,z^q) \}_{q=p^e}$ is cofinal with the direct
system $\{ S/(x^n,y^n,z^n) \}_{n \in \N}$, and so we have
$$
H^3_{(x,y,z)}(S) \cong \varinjlim_{q=p^e} S/(x^q,y^q,z^q)S.
$$
Using this, it is easily seen that
$$
\Ass H^3_{(x,y,z)}(S) \subseteq \bigcup_{q=p^e} \Ass S/(x^q,y^q,z^q)S.
$$
By Theorem \ref{main} $H^3_{(x,y,z)}(S)$ has infinitely many associated prime
ideals, and so $\bigcup_{q}\Ass S/(x^q,y^q,z^q)S$ must be infinite as well.
Since the hypersurface $S$ is F-regular, we have $(x^q,y^q,z^q)^*=(x^q,y^q,z^q)$
for all $q=p^e$. 
\end{proof}

\begin{remark}
In \cite{SS} we actually prove a stronger result: There exists an F-regular
hypersurface $R$ of characteristic $p>0$, with an ideal $\a$, for which the set
$$
\bigcup_{q=p^e} \Ass R/\a^{[q]}
= \bigcup_{q=p^e} \Ass R/\big(\a^{[q]}\big)^*
$$
has infinitely many {\it maximal}\/ elements.
\end{remark}

\end{document}